\documentclass[10pt]{amsart}
\usepackage{latexsym}   
\usepackage{amssymb}    
\usepackage{amsmath}    
\usepackage{amsthm}     
\usepackage[all]{xy}
\usepackage{hyperref}
\newcommand{\pa}{\partial}

\newcommand{\del}{\delta}

\newcommand{\la}{\lambda}\newcommand{\La}{\Lambda}\newcommand{\om}{\omega}

\newcommand{\ti}{\tilde}

\renewcommand{\thefootnote}

\newtheorem{theorem}{Theorem}[section]

\theoremstyle{definition}

\theoremstyle{remark}

\numberwithin{equation}{section}

\title[Bianchi's Hazzidakis transformation for generic higher dimensional quadrics] {
Bianchi's Hazzidakis transformation for generic higher dimensional quadrics}

\author[  Ion I. Dinc\u{a}]{Ion I. Dinc\u{a}}
\address{Department of Applied Mathematics, Faculty of Applied Sciences,
University Politehnica of Bucharest 313 Spl. Independentei 060042
Bucharest, Romania}
 \email{idinca@upb.ro}

\begin{document}


\begin{abstract}
We provide a generalization of Bianchi's Hazzidakis transformation from $2$-dimensional quadrics to generic higher dimensional quadrics.
\end{abstract}

\maketitle

\tableofcontents \pagenumbering{arabic}

\section{Introduction}

For more details on the history of isometric deformations of surfaces and (higher dimensional) quadrics see
Dinc\u{a} \cite{D1}.

All computations are local and assumed to be valid on their open domain of validity without further
details; all functions have the assumed order of differentiability and are assumed to be invertible,
non-zero, etc when required (for all practical purposes we can assume all functions to be analytic).

Consider the complexified Euclidean space
$$(\mathbb{C}^p,<.,.>),\ <x,y>:=x^Ty,\ |x|^2:=x^Tx,\ x,y\in\mathbb{C}^p$$
with standard basis $\{e_j\}_{j=1,...,p},\ e_j^Te_k=\del_{jk}$.

Isotropic (null) vectors are those vectors $v$ of length $0$ ($|v|^2=0$); since most vectors are not
isotropic we shall call a vector simply vector and we shall only emphasize isotropic when the vector
is assumed to be isotropic. The same denomination will apply in other settings: for example we call
quadric a non-degenerate quadric (a quadric projectively equivalent to the complex unit sphere).

A quadric $x\subset\mathbb{C}^{n+1}$ is given by the quadratic equation
$$Q(x):=\begin{bmatrix}x\\1\end{bmatrix}^T\begin{bmatrix}A&B\\B^T&C\end{bmatrix}
\begin{bmatrix}x\\1\end{bmatrix}=x^T(Ax+2B)+C=0,$$
$$A=A^T\in\mathbf{M}_{n+1}(\mathbb{C}),\ B\in\mathbb{C}^{n+1},\ C\in\mathbb{C},\
\begin{vmatrix}A&B\\B^T&C\end{vmatrix}\neq 0.$$

The quadric $x\subset\mathbb{C}^{n+1}$ is degenerate if $\begin{vmatrix}A&B\\B^T&C\end{vmatrix}=0$.

There are many definitions of totally real (sub)spaces of $\mathbb{C}^{n+1}$, some even involving a
hermitian inner product, but all definitions coincide: an $(n+1)$-totally real subspace of
$\mathbb{C}^{n+1}$ is of the form $(R,t)(\mathbb{R}^k\times(i\mathbb{R})^{n+1-k}),\ k=0,...,n+1$,
where $(R,t)\in\mathbf{O}_{n+1}(\mathbb{C})\ltimes\mathbb{C}^{n+1}$. Now a totally real quadric is
simply an $n$-dimensional quadric in an $(n+1)$-totally real subspace of $\mathbb{C}^{n+1}$.

A metric classification of all (totally real) quadrics in $\mathbb{C}^{n+1}$ requires the notion of
symmetric Jordan canonical form of asymmetric complex matrix. The symmetric Jordan blocks
are:
$$J_1:=0=0_{1,1}\in\mathbf{M}_1{\mathbb{C}},\ J_2:=f_1f_1^T\in\mathbf{M}_2{\mathbb{C}},\
J_3:=f_1e_3^T+e_3f_1^T\in\mathbf{M}_3{\mathbb{C}},$$
$$J_4:=f_1\bar f_2^T+f_2f_2^T+\bar f_2f_1^T\in\mathbf{M}_4{\mathbb{C}},\
J_5:=f_1\bar f_2^T+f_2e_5^T+e_5f_2^T+\bar f_2f_1^T\in\mathbf{M}_5{\mathbb{C}},$$
$$J_6:=f_1\bar f_2^T+f_2\bar f_3^T+f_3f_3^T+\bar f_3f_2^T+\bar f_2f_1^T\in\mathbf{M}_6{\mathbb{C}},$$
etc,where $f_j:=\frac{e_{2j-1}-ie_{2j}}{\sqrt{2}}$ are the standard isotropic vectors (at least the
blocks $J_2,\ J_3$ were known to the classical geometers).

Any symmetric complex matrix can be brought via conjugation with a complex rotation to the symmetric
Jordan canonical form, that is a matrix block decomposition with blocks of the form $a_jI_p+J_p$;
totally real quadrics are obtained for eigenvalues $a_j$ of the quadratic part $A$ defining the
quadric being real or coming in complex conjugate pairs $a_j,\ \bar a_j$ with subjacent symmetric
Jordan blocks of same dimension $p$. Just as the usual Jordan block $\sum_{j=1}^pe_je_{j+1}^T$ is
nilpotent with $e_{p+1}$ cyclic vector of order $p$, $J_p$ is nilpotent with $\bar f_1$ cyclic
vector of order $p$, so we can take square roots of symmetric Jordan matrices without isotropic kernels
($\sqrt{aI_p+J_p}:=\sqrt{a}\sum_{j=0}^{p-1}(_j^{\frac{1}{2}})a^{-j}J_p^j,\ a\in\mathbb{C}^*,\
\sqrt{a}:=\sqrt{r}e^{i\theta}$ for $a=re^{2i\theta},\ 0<r,\ -\pi<2\theta\le\pi$), two matrices with
same symmetric Jordan decomposition type (that is $J_p$ is replaced with a polynomial in $J_p$) commute, etc.

The confocal family $\{x_z\}_{z\in\mathbb{C}}$ of a quadric with center $x_0\subset\mathbb{C}^{n+1}$ in
canonical form (depending on as few constants as possible) is given in the projective space
$\mathbb{CP}^{n+1}$ by the equation
$$Q_z(x_z):=\begin{bmatrix}x_z\\1\end{bmatrix}^T(\begin{bmatrix}A&0\\0^T&-1\end{bmatrix}^{-1}-
z\begin{bmatrix}I_{n+1}&0\\0^T&0\end{bmatrix})^{-1}\begin{bmatrix}x_z\\1\end{bmatrix}=0,$$
where $A=A^T\in\mathbf{GL}_{n+1}(\mathbb{C})$ is symmetric Jordan.

From the definition one can see that the family of quadrics confocal to $x_0$ is the adjugate of the
pencil generated by the adjugate of $x_0$ and Cayley's absolute $C(\infty)\subset\mathbb{CP}^n$ in
the hyperplane at infinity; since Cayley's absolute encodes the Euclidean structure of
$\mathbb{C}^{n+1}$ (it is the set invariant under rigid motions and homotheties of
$\mathbb{C}^{n+1}:=\mathbb{CP}^{n+1}\setminus\mathbb{CP}^n$) the mixed {\it metric-projective}
character of the confocal family becomes clear.

With $R_z:=I_{n+1}-zA,\ z\in\mathbb{C}\setminus\mathrm{spec}(A)^{-1}$ the family of quadrics
$\{x_z\}_z$ confocal to the quadric with center $x_0$ is given by $Q_z(x_z)=x_z^TAR_z^{-1}x_z-1=0$.

The Ivory affine transformation for confocal quadrics with center is an affine correspondence between confocal quadrics with center and having
good metric properties: it is given by
$$x_z=\sqrt{R_z}x_0.$$

For the specific computations of isometric deformations of quadrics with center we shall use the convention
$\mathbb{C}^n\subset\mathbb{C}^{n+1}$ with $0$ on the $(n+1)^{\mathrm{th}}$ component; thus for
example we can multiply $(n+1,n+1)$-matrices with $n$-column vectors and similarly one can extend
$(n,n)$-matrices to $(n+1,n+1)$-matrices with zeroes on the last column and row. The converse is
also valid: an $(n+1,n+1)$-matrix with zeroes on the last column and row (or multiplied on the left
with an $n$-row vector and on the right with an $n$-column vector) will be considered as an
$(n,n)$-matrix.

With $V:=\sum_{k=1}^nv^ke_k=[v^1\ ...\ v^n]^T$ consider the complex unit sphere
$$X=X(v^1,...,v^n)=\frac{2V+(|V|^2-1)e_{n+1}}{|V|^2+1}.$$
We have
$$dX=2\frac{dV+V^TdV(e_{n+1}-X)}{|V|^2+1},\ \mathrm{so}\ |dX|^2=\frac{4|dV|^2}{(|V|^2+1)^2}.$$

With $A=A^T\in\mathbf{GL}_{n+1}(\mathbb{C})$\ symmetric Jordan we have the quadric with center $x_0=(\sqrt{A})^{-1}X$.

Consider the non-trivial isometric deformation $x\subset\mathbb{C}^{2n-1}$ of the quadric with center $x_0\subset\mathbb{C}^{n+1}$ with common conjugate system $(u^1,...,u^n)$.

Because $(v^1,...,v^n)$ are isothermal-conjugate and $(u^1,...,u^n)$ are conjugate on $x_0$, the
Jacobian $\frac{\pa(v^1,...,v^n)}{\pa(u^1,...,u^n)}$ has orthogonal columns, so with
$$\la_j:=|\pa_{u^j}V|\neq 0,\ \La:=[\la_1\ ...\ \la_n]^T,\ \del:=\mathrm{diag}[du^1\ ...\ du^n]$$
we have $dV=R\del\La,\ R\subset\mathbf{O}_n(\mathbb{C})$.

With $I_{1,n}:=I_{n+1}-e_{n+1}e_{n+1}^T$ we have the next theorem from Dinc\u{a} \cite{D1} on non-trivial isometric deformations of quadrics with center

\begin{theorem}\label{th:th5}
 Non-trivial isometric deformations $x\subset\mathbb{C}^{2n-1}$ of the quadric with center $x_0\subset\mathbb{C}^{n+1}$ with the isothermal-conjugate system $(v^1,...,v^n)$ and common conjugate system $(u^1,...,u^n)$ are in
bijective correspondence with solutions of the differential system (\ref{eq:domqc}) in involution (that is no further conditions appear if one imposes $d$ conditions and one uses the
equations of the system itself).
\end{theorem}

\begin{eqnarray}\label{eq:domqc}
d\om-\om\wedge\om=-4\del R^T(I_{1,n}+Ve_{n+1}^T)A(I_{1,n}+e_{n+1}V^T)R\wedge\del,\
\om\wedge\del-\del\wedge R^TdR=0,\nonumber\\
dV=R\del\La,\ d\La=\om\La-2\del R^T(I_{1,n}+Ve_{n+1}^T)A(2V+(|V|^2-1)e_{n+1}),\nonumber\\
\La^T\La=-(2V+(|V|^2-1)e_{n+1})^TA(2V+(|V|^2-1)e_{n+1}),\nonumber\\
\om:=\sum_{j=1}^n(e_je_j^TR^T\pa_{u^j}R\del+\del R^T\pa_{u^j}Re_je_j^T),\
R\subset\mathbf{O}_n(\mathbb{C}).
\end{eqnarray}

We have the next theorem from Dinc\u{a} \cite{D1} on the B\"{a}cklund transformation for higher dimensional quadrics

\begin{theorem}\label{th:th1}
{\bf I (existence and inversion of the B\"{a}cklund transformation for quadrics and the isometric
correspondence provided by the Ivory affine transformation)}

Any non-trivial isometric deformation $x^0\subset\mathbb{C}^{2n-1}$ of an $n$-dimensional sub-manifold
$x_0^0\subset x_0$ ($x_0\subset\mathbb{C}^{n+1}\subset\mathbb{C}^{2n-1}$ being a quadric) appears as
a focal sub-manifold of an $(\frac{n(n-1)}{2}+1)$-dimensional family of Weingarten congruences, whose
other focal sub-manifolds $x^1=B_z(x^0)$ are isometric, via the Ivory affine transformation between
confocal quadrics, to sub-manifolds $x_0^1$ in the same quadric $x_0$. The determination of these
sub-manifolds requires the integration of a family of Ricatti equations depending on the parameter
$z$ (we ignore for simplicity the dependence on the initial value data in the notation $B_z$).
Moreover, if we compose the inverse of the rigid motion provided by the Ivory affine transformation $(R_0^1,t_0^1)$ with the rolling of $x_0^0$ on $x^0$, then we obtain the rolling of $x_0^1$
on $x^1$ and $x^0$ reveals itself as a $B_z$ transformation of $x^1$.
\end{theorem}

Recall that we have the Ivory affine transformation $x_z=\sqrt{R_z}(\sqrt{A})^{-1}X,\
X=\frac{2V+(|V|^2-1)e_{n+1}}{|V|^2+1}.$

Consider two points $x_0^0,\ x_0^1\in x_0$ such that $x_0^0,\ x_z^1$ are in the symmetric tangency configuration

\begin{eqnarray}\label{eq:tcqc}
(x_z^1-x_0^0)^TN_0^0=0\Leftrightarrow X_1^T\sqrt{R_z}X_0=1\Leftrightarrow x_z^1=x_0^0+
\frac{(|V_0|^2+1)^2}{4}\sum_{k=1}^nX_1^T\sqrt{R_z}\pa_{v_0^k}X_0\pa_{v_0^k}x_0^0\Leftrightarrow
\nonumber\\
(2V_0+(|V_0|^2-1)e_{n+1})^T\sqrt{R_z}(2V_1+(|V_1|^2-1)e_{n+1})=(|V_0|^2+1)(|V_1|^2+1).
\end{eqnarray}
Thus among the $2n$ functionally independent variables $\{v_0^j,v_1^j\}_{j=1,...,n}$ a quadratic
functional relationship is established and only $2n-1$ among them remain functionally independent

\begin{eqnarray}\label{eq:X1Rz}
X_1^T\sqrt{R_z}dX_0=-X_0^T\sqrt{R_z}dX_1.
\end{eqnarray}
Given a non-trivial isometric deformation $x^0\subset\mathbb{C}^{2n-1}$ of $x_0^0$ (that is $|dx^0|^2=
|dx_0^0|^2$) with orthonormal normal frame $N^0:=[N_{n+1}^0\ ...\ N_{2n-1}^0]$ consider the
$n$-dimensional sub-manifold

\begin{eqnarray}\label{eq:x1x0qc}
x^1=x^0+\frac{(|V_0|^2+1)^2}{4}\sum_{k=1}^nX_1^T\sqrt{R_z}\pa_{v_0^k}X_0\pa_{v_0^k}x^0\subset
\mathbb{C}^{2n-1}
\end{eqnarray}
(that is we restrict $\{v_1^j\}_{j=1,...,n}$ to depend only on the functionally independent
$\{v_0^j\}_{j=1,...,n}$ and constants in a manner that will subsequently become clear when we shall
impose the isometric correspondence provided by the Ivory affine transformation).

We take advantage now of the conjugate system $(u^1,...,u^n)$ common to $x_0^0,\ x^0$ and of the
non-degenerate joined second fundamental forms property; according to the principle of symmetry
$0\leftrightarrow 1$ we would like $(u^1,...,u^n)$ to be conjugate system to both $x^1$ and $x_0^1$
and also that the non-degenerate joined second fundamental forms property holds for $x_0^1,\ x^1$.
We obtain $dV_1=R_1\del\La_1$,
$$\La_1:=\frac{(|V_1|^1+1)(|V_0|^2+1)}{2\sqrt{z}}R_1^T[X_1^T\sqrt{R_z}\pa_{v_0^1}X_0\ ...
\ X_1^T\sqrt{R_z}\pa_{v_0^n}X_0]^T$$
(note that the prime integral property $|\La_1|^2=-(2V_1+(|V_1|^2-1)e_{n+1})^TA(2V_1+(|V_1|^2-1)e_{n+1})$ is satisfied).

Thus

\begin{eqnarray}\label{eq:zR0La1}
\sqrt{z}R_0\La_1=(I_{1,n}+V_0e_{n+1}^T)\sqrt{R_z}(2V_1+(|V_1|^2-1)e_{n+1})-V_0(|V_1|^2+1),\nonumber\\
-\sqrt{z}R_1\La_0=(I_{1,n}+V_1e_{n+1}^T)\sqrt{R_z}(2V_0+(|V_0|^2-1)e_{n+1})-V_1(|V_0|^2+1)
\end{eqnarray}
(the symmetry $(0,\sqrt{z})\leftrightarrow(1,-\sqrt{z})$ follows from (\ref{eq:X1Rz})). From the
second equation of (\ref{eq:zR0La1}) we get $V_1$ as an algebraic expression of $R_1,\ V_0,\ \La_0$;
replacing this into the first relation of (\ref{eq:zR0La1}) we get $\La_1$ as an algebraic
expression of $R_0,\ R_1,\ V_0,\ \La_0$:

\begin{eqnarray}\label{eq:V1La1qc}
V_1=-\frac{\sqrt{z}R_1\La_0+I_{1,n}\sqrt{R_z}(2V_0+(|V_0|^2-1)e_{n+1})}
{e_{n+1}^T\sqrt{R_z}(2V_0+(|V_0|^2-1)e_{n+1})-|V_0|^2-1},\nonumber\\
\La_1=2R_0^T\frac{(I_{1,n}+V_0e_{n+1}^T)[\sqrt{z}A(2V_0+(|V_0|^2-1)e_{n+1})
-\sqrt{R_z}(I_{1,n}+e_{n+1}V_1^T)R_1\La_0]+V_0V_1^TR_1\La_0}
{e_{n+1}^T\sqrt{R_z}(2V_0+(|V_0|^2-1)e_{n+1})-|V_0|^2-1}.\nonumber\\
\end{eqnarray}

Differentiating the second relation of (\ref{eq:zR0La1}) we get
$$-\sqrt{z}dR_1\La_0-\sqrt{z}R_1[\om_0\La_0-2\del R_0^T(I_{1,n}+V_0e_{n+1}^T)A(2V_0+(|V_0|^2
-1)e_{n+1})]-$$
$$2I_{1,n}\sqrt{R_z}(I_{1,n}+e_{n+1}V_0^T)R_0\del\La_0=
[e_{n+1}^T\sqrt{R_z}(2V_0+(|V_0|^2-1)e_{n+1})-|V_0|^2-1]R_1\del\La_1+$$
$$2V_1[e_{n+1}^T\sqrt{R_z}(I_{1,n}+e_{n+1}V_0^T)-V_0^T]R_0\del\La_0,$$
or equivalently the Ricatti equation

\begin{eqnarray}\label{eq:dR1qc}
-dR_1=R_1\om_0+2I_{1,n}\frac{\sqrt{R_z}}{\sqrt{z}}(I_{1,n}+e_{n+1}V_0^T)R_0\del
-2R_1\del R_0^T(I_{1,n}+V_0e_{n+1}^T)\frac{\sqrt{R_z}}{\sqrt{z}}R_1\nonumber\\
+2R_1\del R_0^T\frac{[(I_{1,n}+V_0e_{n+1}^T)\sqrt{R_z}e_{n+1}-V_0][\La_0^T+(2V_0^T+
(|V_0|^2-1)e_{n+1}^T)\frac{\sqrt{R_z}}{\sqrt{z}}R_1]}{e_{n+1}^T\sqrt{R_z}(2V_0+(|V_0|^2-1)e_{n+1})
-|V_0|^2-1}\nonumber\\
-2\frac{[R_1\La_0+I_{1,n}\frac{\sqrt{R_z}}{\sqrt{z}}(2V_0+(|V_0|^2-1)e_{n+1})]
[e_{n+1}^T\sqrt{R_z}(I_{1,n}+e_{n+1}V_0^T)-V_0^T]}{e_{n+1}^T\sqrt{R_z}(2V_0+(|V_0|^2-1)e_{n+1})
-|V_0|^2-1}R_0\del
\end{eqnarray}
in $R_1$.

We have now from Dinc\u{a} \cite{D1} the theorem partially corresponding to {\bf Theorem} \ref{th:th1} for quadrics with center:

\begin{theorem}\label{th:th7}
Given $(V_0,\La_0,R_0)$ solution of (\ref{eq:domqc}) that produces a non-trivial isometric deformation
$x^0\subset\mathbb{C}^{2n-1}$ of $x_0^0\subset x_0$, $x_0$ being a quadric with center the Ricatti
equation (\ref{eq:dR1qc}) is completely integrable. If $R_1\subset\mathbf{M}_n(\mathbb{C})$ is a
solution of (\ref{eq:dR1qc}), then it remains orthogonal if initially it was orthogonal and in this
case together with $V_1,\ \La_1$ given by (\ref{eq:V1La1qc}) it is another solution of
(\ref{eq:domqc}) (thus producing an $(\frac{n(n-1)}{2}+1)$-dimensional family of leaves non-trivial isometric
deformations $x^1\subset\mathbb{C}^{2n-1}$ of $x_0^1$) and we have the symmetry $(0,\sqrt{z})
\leftrightarrow(1,-\sqrt{z})$.
\end{theorem}

Note that with
$$M_0(z):=I_{1,n}\frac{\sqrt{R_z}}{\sqrt{z}}(I_{1,n}+e_{n+1}V_0^T),\
P_0(z):=I_{1,n}\frac{\sqrt{R_z}}{\sqrt{z}}(2V_0+(|V_0|^2-1)e_{n+1}),$$
$$W_0(z):=(I_{1,n}+V_0e_{n+1}^T)\sqrt{R_z}e_{n+1}-V_0,\
U_0(z):=e_{n+1}^T\sqrt{R_z}(2V_0+(|V_0|^2-1)e_{n+1})-|V_0|^2-1$$
we have

$$dP_0(z)=2M_0(z)dV_0,\ dU_0(z)=2W_0(z)^TdV_0,\ dM_0(z)=C(z)dV_0^T,\ dW_0(z)=c(z)dV_0$$
and (\ref{eq:dR1qc}) can be written as

\begin{eqnarray}\label{eq:dR1qcm}
-dR_1=R_1\om_0+2M_0(z)R_0\del-2R_1\del R_0^TM_0(z)^TR_1
+\frac{2}{U_0(z)}R_1\del R_0^TW_0(z)(\La_0^T+P_0(z)^TR_1)\nonumber\\
-\frac{2}{U_0(z)}(R_1\La_0+P_0(z))W_0(z)^TR_0\del.
\end{eqnarray}

We have now the main theorem of this paper:

\begin{theorem}\label{th:th8}
{\bf IV (Hazzidakis transformation)}
If an $n$-dimensional sub-manifold $x^0\subset\mathbb{C}^{2n-1}$ is a non-trivial isometric deformation of a sub-manifold $x_0^0\subseteq x_0$, $x_0\subset\mathbb{C}^{n+1}$ being a generic quadric and the homography $H\in\mathbf{PGL}_{n+1}(\mathbb{C})$ takes the confocal family $x_z$ to another confocal family
$\ti x_{\ti z},\ \ti z=\ti z(z),\ \ti z(0)=0$ of a generic quadric $\ti x_0\subset\mathbb{C}^{n+1}$, then one infinitesimally knows a sub-manifold $\ti x^0=H(x^0)\subset\mathbb{C}^{2n-1}$ (that is one knows the first and second fundamental forms)
called the Hazzidakis transform of $x^0$ and non-trivial isometric deformation of a sub-manifold $\ti x_0^0\subseteq\ti x_0$. Moreover
the Hazzidakis transformation commutes with the B\"{a}cklund transformation ($H\circ B_z=B_{\ti z}\circ H$) and the $B_{\ti z}(\ti x^0)$ transforms can be algebraically recovered from the knowledge of $\ti x^0$ and $B_z(x^0)$.
\end{theorem}

In the next section we provide the proof of the main theorem of this paper.

\section{Proof of the Hazzidakis transformation}

A generic quadric $x_0\subset\mathbb{C}^{n+1}$ is a quadric with center given by $A=\mathrm{diag}[a_1^{-1}\ a_2^{-1}\ ...\ a_{n+1}^{-1}]$ with $a_j\in\mathbb{C}\setminus\{0\},\ j=1,...,n+1$ distinct.

For a generic higher dimensional quadric with $A=\mathrm{diag}[a_1^{-1}\ a_2^{-1}\ ...\ a_{n+1}^{-1}]$ (\ref{eq:domqc}) becomes

\begin{eqnarray}\label{eq:domqc1}
d\om-\om\wedge\om=-4\del R^T(I_{1,n}AI_{1,n}+a_{n+1}^{-1}VV^T)R\wedge\del,\
\om\wedge\del-\del\wedge R^TdR=0,\nonumber\\
dV=R\del\La,\ d\La=\om\La-2\del R^T(2I_{1,n}AI_{1,n}V+a_{n+1}^{-1}(|V|^2-1)V),\nonumber\\
\La^T\La=-(4V^TI_{1,n}AI_{1,n}V+a_{n+1}^{-1}(|V|^2-1)^2),\nonumber\\
\om:=\sum_{j=1}^n(e_je_j^TR^T\pa_{u^j}R\del+\del R^T\pa_{u^j}Re_je_j^T),\
R\subset\mathbf{O}_n(\mathbb{C})
\end{eqnarray}
and (\ref{eq:zR0La1}) becomes

\begin{eqnarray}\label{eq:zR0La11}
\sqrt{z}R_0\La_1=2I_{1,n}\sqrt{R_z}I_{1,n}V_1+\sqrt{1-za_{n+1}^{-1}}(|V_1|^2-1)V_0-V_0(|V_1|^2+1),\nonumber\\
-\sqrt{z}R_1\La_0=2I_{1,n}\sqrt{R_z}I_{1,n}V_0+\sqrt{1-za_{n+1}^{-1}}(|V_0|^2-1)V_1-V_1(|V_0|^2+1).
\end{eqnarray}

In the following tilde quantities correspond to non-tilde quantities and conversely.

The generic quadric $\ti x_0\subset\mathbb{C}^{n+1}$ is a quadric with center given by $\ti A=\mathrm{diag}[\ti a_1^{-1}\ \ti a_2^{-1}\ ...\ \ti a_{n+1}^{-1}]$, where $\ti a_{n+1}=a_{n+1}^{-1},\ \ti a_j=\ti a_{n+1}+(a_j-a_{n+1})^{-1},\ j=1,...,n,\ \ti z=\ti a_{n+1}+(z-a_{n+1})^{-1}$.

The homography $H$ which takes the confocal family $x_z$ of $x_0$ to the confocal family $\ti x_{\ti z}$ of $\ti x_0$ is given by $\mathbb{C}^{n+1}\ni x\mapsto H(x)=\frac{H'x+e_{n+1}}{e_{n+1}^Tx},\ H'=\sum_{j=1}^n\frac{1}{\sqrt{a_{n+1}-a_j}}e_je_j^T,\ \ti H'=\sum_{j=1}^n\sqrt{a_{n+1}-a_j}e_je_j^T$.

We have $x_0=(\sqrt{A})^{-1}\frac{2V_0+(|V_0|^2-1)e_{n+1}}{|V_0|^2+1},\ \ti x_0=(\sqrt{\ti A})^{-1}\frac{2\ti V_0+(|\ti V_0|^2-1)e_{n+1}}{|\ti V_0|^2+1}$, so $\ti V_0=iV_0$ (in fact we have $\ti v_0^j=\pm iv_0^j,\ j=1,...,n$, but we
assume simplifications of the form $\sqrt{a}\sqrt{b}\simeq\sqrt{ab},\ a,b\in\mathbb{C}$ since the possible signs are accounted by symmetries in the principal hyper-planes and disappear at the level of the linear element for isometric deformations).

Note that $I_{1,n}\ti AI_{1,n}=a_{n+1}I_{1,n}-a_{n+1}^2I_{1,n}AI_{1,n}$.

If we choose $\ti u^j:=ia_{n+1}^{-1}u^j,\ j=1,...,n$, then $\pa_{\ti u^j}=\frac{1}{ia_{n+1}^{-1}}\pa_{u^j},\
d\ti u^j=ia_{n+1}^{-1}du^j$, so $\ti\del=ia_{n+1}^{-1}\del$.

Choose $\ti R_0:=R_0$; from $d\ti V_0=\ti R_0\ti\del\ti\La_0$ we get $\ti\La_0=a_{n+1}\La_0$. From the definition we have $\ti\om_0=\om_0$ and (\ref{eq:domqc1}) is satisfied with $\ti R_0,\ti\om_0,\ti\del,\ti V_0,\ti\La_0.\ti A$
replacing respectively $R,\om,\del,V,\La,A$ if (\ref{eq:domqc1}) is satisfied with $R_0,\om_0,\del,V_0,\La_0.A$
replacing respectively $R,\om,\del,V,\La,A$. Thus we know infinitesimally an isometric deformation $\ti x^0$ of
an $n$-dimensional sub-manifold $\ti x_0^0\subseteq\ti x_0$ (that is we know the first and second fundamental forms). 

Assume now that we have a B\"{a}cklund transform $x^1=B_z(x^0)$ of $x^0$; to get a B\"{a}cklund transform $\ti x^1=B_{\ti z}(\ti x^0)$ (that is such that $\ti x^1$ is a Hazzidakis transform of $x^1$) that can be algebraically recovered from the knowledge of
$x^1$ and $\ti x^0$ we need to show that (\ref{eq:zR0La11}) is satisfied with tilde quantities replacing non-tilde quantities if it is satisfied initially for non-tilde quantities(in this case the equation (\ref{eq:dR1qcm}) obtained by applying $d$ to the second equation of (\ref{eq:zR0La11}) needs not to be checked for tilde quantities replacing non-tilde quantities).

We have $\ti z=-\frac{z}{a_{n+1}^2}(1-za_{n+1}^{-1})^{-1},\ 1-\ti z\ti a_{n+1}^{-1}=(1-za_{n+1}^{-1})^{-1}$ and 

$I_{1,n}\ti R_{\ti z}I_{1,n}=(1-za_{n+1}^{-1})^{-1}I_{1,n}R_zI_{1,n}$, so (\ref{eq:zR0La11}) is satisfied with tilde quantities replacing non-tilde quantities if it is satisfied initially for non-tilde quantities (the signs on the left hand of (\ref{eq:zR0La11}) may need to exchange).

\end{document}